\documentclass[12pt]{article}
\usepackage[mathscr]{eucal}
\usepackage{amssymb,amsmath}
\font\cmc=cmcsc10  scaled \magstep2
\newcommand\de{\delta}

\newcommand\va{\varphi}
\newcommand{\ve}{\varepsilon}

\newcommand\la{\lambda}

\newcommand\n{\noindent}

\newcommand\vk{\vskip}

\newcommand\al{\alpha}

\newcommand\Ho{\text{Hol}}

\newcommand{\Hol}{\text{Hol}}

\newcommand\nhd{neigborhood}
\newcommand\De{\Delta}

\newcommand\psh{plurisubharmonic}
\newcommand\no{\noindent}

\newcommand\cen{\centerline}
\newcommand\Om{\Omega}
\newcommand\ovl{\overline}

\newtheorem{Proof.}{\it Proof.}

\begin{document}
\vbox to .5truecm{}
\large
\begin{center}
\cmc Vitali properties of Banach analytic manifolds
\end{center}
\centerline {Nguyen Van Khue, Nguyen Quang Dieu and Nguyen Van Khiem}
\vk.5cm
{\normalsize
\no{\bf Abstract.} We discuss possible generalizations of Vitali convergence theorem when the source and the target are Banach analytic manifolds. These results are then applied to study behavior of holomorphic mappings between Banach analytic manifolds. Explicit examples of manifolds having Vitali properties are also provided.
\vk.3cm
\no {\bf Keywords} Vitali theorem, hyperbolic manifold, taut property, Banach analytic manifold

\no {\bf Mathematics Subject Classification 2010:} 32A10, 32A07, 32C25, 32Q45, 32U05
}
\vk.5cm

\cen{\textbf {I. Introduction}}
\vk.5cm
\no
The classical Vitali theorem states that a sequence $\{f_k\}$ of holomorphic functions defined on a domain $D$ in $\mathbb C$ is uniformly convergent on compact sets if it is locally uniformly bounded and
if it converges {\it pointwise} only on some set having an accumulation point in $D$.
There are two ingredients in the proof. Firstly, by Montel's theorem, the sequence $\{f_k\}$ is relatively compact in the compact open topology and secondly, using the uniqueness property of holomorphic functions, we conclude that two accumulation points of the sequence $\{f_k\}$ must coincide on $D$. Observe that it is rather straightforward to generalize Vitali's theorem to (scalar valued) holomorphic functions of several variable.
For vector-valued holomorphic functions, Montel's theorem is not valid, therefore it does not seem easy to find an analogue of Vitali theorem in this more general setting.
Nevertheless,  by making use the notion of weak holomorphicity together with some elementary but quite ingenious arguments, Arendt and Nikolski provide in [1] a correct generalization of Vitali theorem for holomorphic functions defined on domains in $\mathbb C$ with values in Banach spaces.
The aim of this paper is to explore possible versions of Vitali theorems in a general setting where the source and the target spaces are assumed to be Banach analytic manifolds.

Now, we will shortly review basic notions that pertaining to our work.
By a Banach analytic manifold we mean a {\it connected} topological space in which each point has a neighborhood homeomorphic to an open set in a Banach space such that the transition maps are holomorphic between open sets of Banach spaces.
Thus, Banach analytic manifolds encompass two objects of different character: (finite dimensional) complex manifolds and
(infinite dimensional) Banach space.

Roughly speaking, we say a Banach analytic manifold $X$ has the Vitali property if for every (connected) Banach analytic manifold $A$ and every sequence $\{f_k\}$ of holomorphic mappings from $A$ into $X$ that converges only pointwise on a "sufficiently large" subset of $A$ must converge uniformly on compact sets of $A.$
For clarity of the exposition, we introduce Banach analytic manifolds with weak Vitali property (WVP) and strong Vitali property (SVP) depending on the nature of the set where pointwise convergence of $\{f_k\}$ occurs.
Even though, we do not know if the two properties are really different, there are certainly some advantages in studying them. We now briefly outline the content of the paper. The first part concentrates on theoretic properties of manifolds having Vitali properties and their applications to study behavior of holomorphic mappings between Banach analytic manifolds.
Our first main results is Theorem 3.1 which says that every Banach analytic manifold having WVP must be (Kobayashi) hyperbolic. This result brings in naturally hyperbolic Banach analytic manifolds into our study.
In the opposite direction, we show in Theorem 3.3 that every complete hyperbolic Banach analytic
manifold has SVP. The proof relies strongly on a vector valued version of Vitali's theorem which is inspired from the mentioned above work of Arendt and Nikolski in [1].
We also relate our Vitali properties with some sorts of {\it taut}  property of Banach analytic manifolds.
Recall that the classical taut property (see [7], p. 239) is defined for (finite dimensional) complex manifolds  and it reflect the behavior of sequences of holomorphic mappings from the unit disk $\De \subset \mathbb C$
into the complex manifold under consideration.

Our Vitali properties serves as convenient tools to check tautness of complex manifolds and Banach analytic manifolds. This fact is reflected in Theorem 3.10 which says that every sequence of holomorphic maps from a connected separable Banach analytic manifold $A$ into a Banach analytic manifold $X$ having WVP must contain a subsequence which is either convergent or compactly divergent on an open {\it dense} subset of $A$.
Under the stronger assumption that the target manifold $X$ has SVP and the source manifold is just the unit disk $\De$, we show in Theorem 3.11 that the compactly divergence phenomenon {\it may} only occur outside a discrete subset of $\De.$
We should mention that, in the literature, there are some attempts to generalize the classical taut property for Banach analytic manifolds (see [4], [5] and [6]). It should be, however, noted that our proofs, unlike those in [4], [5], [6], are quite constructive, in the sense that we avoid to use Zorn's lemma.

The second part of the work focus on finding explicit classes of manifolds $X$ having WVP and SVP.
The key idea is to impose on the existence of certain (non-constant) negative \psh\ on $X$, and under certain additional
assumptions we get Vitali property of the whole space $X$ if each sublevel set determined by $\va$ has this property.
This principle is carried out in Theorem 4.1 (for WVP of Banach analytic manifolds) and Theorem 4.4 (for SVP of open subsets of Banach space). Furthermore, we also give in the last two results somewhat complete characterizations
for Vitali properties of Hartogs domains (over Banach analytic manifolds) and balanced domains in Banach space.
The paper ends up by giving a list of open questions that are connected to our work.
\vskip0,4cm
\noindent
{\bf Acknowledgments.} This work is supported by the grant 101.02-2016.07 from the NAFOSTED program.
\vk1cm
\cen{\textbf {II. Basic notions and notation}}
\vk.2cm
\no
We introduce below certain notion that are needed for formulating Vitali properties.
\vk.3cm
\no{\bf Notation.} (a) Let $A$ be a Banach analytic manifold and $S$ be a subset of $A$. We let

\no
$S^u:=\Big\{z\in A \cap\ \ovl S: \forall$ connected neighborhood $U$ of $z$ and every holomorphic function $f: U\longrightarrow \mathbb C, f\big|_{ U\cap S}=0 \Rightarrow f\big|_ U=0  \Big\}$.

\noindent
(b) Given Banach analytic manifolds $A$ and $X$. By $\Ho(A, X)$ we mean the linear space of holomorphic mappings from $A$ into $X$. We then equip $\Hol (A,X)$ with the compact-open topology.
According to a result of Palais in [8], a Banach analytic manifold is metrizable if and only if it is paracompact.
So, in the case where $A$ and $X$ are both paracompact, the compact-open topology on $\Hol(A,X)$ is equivalent to the topology of locally uniformly convergence.

\no
(c) Let $A, X$ be Banach analytic manifolds and $\{f_k\}$ be a sequence in $\Hol (A,X)$. We denote by $Z_{\{f_k\}}$ the set of points $\la \in A$ such that $\{f_k (\la)\}$ is convergent.

\noindent
{\bf Remark.} It is easy to check that $S^u$ is a closed subset of $\bar S.$ Moreover, $S \setminus S^u$ is locally contained in an analytic hypersurface i.e., for every $a \in S \setminus S^u$, there exists a connected \nhd\ $U$ of $a$ and a holomorphic function $g$ on $U$ such that $g\not\equiv 0, g|_{S \cap U} \equiv 0.$
\vk.3cm
\no
Now we come the central notions of this paper.
\vk0.2cm
\no\textbf{Definition 2.1}. Let $X$ be a Banach analytic manifold. We say that:

\no
(a) $X$ has the strong Vitali property (SVP for short ) if for every (connected) Banach analytic manifold $A$ and every
sequence $\{f_k\}_{k\geq 1}\subset \text{Hol}(A,X)$ such that $Z^u_{\{f_k\}}\ne \emptyset$ we have $\{f_k\}_{k\geq 1}$ is convergent in $\text{Hol}(A,X)$.

\no
(b) $X$ has the weak Vitali property (WVP for short) if for every (connected) Banach analytic manifold $A$ and every
sequence $\{f_k\}_{k\geq 1}\subset \text{Hol}(A,X)$ such that $Z_{\{f_k\}} \cap Z^u_{\{f_k\}} \ne \emptyset$ we have
$\{f_k\}_{k\geq 1}$ is convergent in $\text{Hol}(A,X)$.

\no
(c) In the particular case where the above properties are true  for $A=\De,$ we say that $X$ has $1-$WVP and $1-$SVP respectively.
\vk.3cm
\no{\bf Remarks 2.2.} (a)
In the special case where $A=\De, Z^u_{\{f_k\}}$ is exactly the set of accumulation points of $Z_{\{f_k\}}.$

\no
(b) We construct a sequence of polynomials $\{p_k\}$ on $\mathbb C$ such that $Z^u_{p_k} \ne \emptyset$ but $Z_{\{p_k\}} \cap Z^u_{\{p_k\}}=\emptyset.$
For $k \ge 1$, we let $U_k$ and $V_k$ be disks in $\mathbb C$ with disjoints closures such that $0 \in U_k$ and $[1/k, 1] \subset V_k.$
By Runge's approximation theorem we can find a polynomial $p_k$ on $\mathbb C$ such that
$$\vert p_k (z)\vert >1/2, \ \forall z \in U_k; \Vert p_k\Vert_{V_k} <1/k.$$
It is then clear that $0 \in Z^u_{\{p_k\}}$ but $0 \not\in Z_{\{p_k\}}.$

\no
(c) In spite of the above example, we will show in Corollary 3.13 and Theorem 3.14 that in the categories of complex manifolds (resp. bounded domains in Banach spaces), the two notions 1-WVP and 1-SVP (resp. WVP and SVP)
are equivalent.
Unfortunately, we do not even know if there exists a (unbounded) domain in a (infinite dimensional) Banach space  having the $1-$WVP but does not have $1-$SVP.
\vk.3cm
\no
The main technical tool in our paper is the Kobayashi pseudo-distance defined on a Banach analytic manifold $X$.
Analogously as in the case where $X$ is a finite dimensional complex manifold (see [7], p.50) or a Banach space (see [3], p. 81), the pseudo-distance $\kappa_X (p, q)$ is defined to be the infimum of the length of all holomorphic chains joining
$p,q \in X.$ More precisely, by a holomorphic chain from
$p$ to $q$ we mean a chain of points $p=p_0, p_1, \cdots, p_k=q$ of $X,$ pairs of points $a_1, b_1, \cdots, a_k, b_k$ of
$\De$ and holomorphic maps $f_1, \cdots, f_k \in \Hol(\De, X)$ such that
$$f_i (a_i)=p_{i-1}, f_i (b_i)=p_i, \ 1 \le i \le k.$$
Denote this chain by $\alpha$, then the length of $\alpha$ is defined to be
$$l(\al):=\rho_\De (a_1, b_1)+\cdots \rho_\De (a_k, b_k),$$
where $\rho_\De$ is the Poincare distance on $\De$. Then the Kobayashi pseudo-distance between $p$ and $q$ is defined by
$$\kappa_X (p, q):=\inf_\alpha l(\alpha),$$
where $\alpha$ is taken over all holomorphic chains connecting $p$ and $q.$

By the same proof as in the case of complex manifolds (see Proposition 3.1.7 in [7]), we can show that $\kappa_X$ is decreasing under holomorphic maps i.e, if $f: X \to Y$ is a holomorphic mapping between Banach analytic manifolds $X, Y$ then
$$\kappa_Y (f(p), q(q)) \le \kappa_X (p, q), \  \forall p, q \in X.$$
Moreover, $\kappa_X$ is the largest pseudo-distance on $X$ having this property.

Then we say that $X$ is {\it hyperbolic} if $\kappa_X$ is a distance and defines the topology of $X$.
Notice that, in contrast to the case where $X$ is finite dimensional,  $\kappa_X$ may be a distance without defining the topology of $X$ even in the case where $X$ is a domain in a Banach space (see [3] p. 93).
Furthermore, $X$ is said to be {\it complete hyperbolic} if every $\kappa_X-$ Cauchy sequence in $X$ is convergent.
By Proposition 6.9 in [2] (see also Proposition 3.6 in [6]) we know that every bounded convex domain $\Om$ in a Banach space is complete hyperbolic.
Hence, all open subsets of $\Om$ are hyperbolic. In particular, each bounded open subset of a Banach space is hyperbolic.
\vk0.2cm
\no
We recall the notion of normality of holomorphic mappings between Banach analytic manifolds when the target space is of finite dimension. This property will be relevant to our Vitali properties in the category of complex manifolds (see Theorem 3.14 in the next section).
\vskip0,2cm
\no
{\bf Definition 2.3.} {\it  Let $A$ be a connected Banach analytic manifolds and $X$ be a complex manifolds. We say that
$\Hol(A,X)$ is normal if for every
sequence $\{f_k\}_{k\geq 1}\subset \text{Hol}(A,X)$ there is a subsequence $\{f_{k_j}\}$ having one of the following properties:

\no
(a) $\{f_{k_j}\}$ is convergent in $\text{Hol}(A,X).$

\no
(b)$\{f_{k_j}\}$ is compactly divergent i.e.,
for every compact subsets $K\subset A$ and $L\subset X$, there exists $j_0$ such that $f_{k_j}(K)\cap L=\emptyset$ for all $j \ge j_0$.}
\vk0.2cm
\no
We will see that the notion of normality does not generalize in the expected fashion when $X$ is a general (infinite dimensional) Banach analytic manifolds. See the remark following Theorem 3.10.

The final ingredient needed in our work is the concept of plurisubharmonic functions on Banach analytic manifolds.
More precisely, we say that $\va: A \to [-\infty, \infty),$ where $A$ is a Banach analytic manifold, is plurisubharmonic if for every $a \in A$, there exists a \nhd\ $U$ of $a$ such that $U$ is isomorphic to a ball $B$ in a Banach space $E$ and that $u,$ regarded as a function on $B,$ is plurisubharmonic in the classical sense (see [3], p.62) i.e, $u$ is upper semicontinuous and the restriction of $u$ on the intersection of $B$ with each complex line in $E$ is subharmonic. Notice that we allow the function $u \equiv -\infty$ to be plurisubharmonic.
We will frequently appeal to the following {\it maximum principle:} Let $A$ be a connected Banach analytic manifold and $\va$ be a \psh\ on $A.$ Suppose that there exists $x_0 \in A$ such that $\va (x_0)=\max_A \va$. Then
$\va|_A \equiv \va(x_0).$

Throughout this paper, for $r>0$, we will write $\De (0,r)$ for the disk in $\mathbb C$ with center $0$ and radius $r.$
\vk1cm
\no
\cen{\textbf {III. Main Results}}
\vk.3cm
\no
Our first result states, in spirit, that hyperbolicity of the target manifold is the right substitute for uniform boundedness assumption given in the classical Vitali theorem.
\vk0.2cm
\no
{\bf Theorem 3.1. } {\it Every Banach analytic manifold $X$ having the $1$-WVP is hyperbolic.}
\vk0,2cm
\noindent
The proof relies heavily on the following lemma which is a slight modification of a result of Kiernan (see Lemma 5.1.4 in [7]).
\vk0,2cm
\noindent
{\bf Lemma 3.2.} {\it Let $Y$ be a Banach analytic manifold and $x \in Y.$ Let $U, V, W$ be open subsets of $Y$ such that
$x \in V \subset \overline{V} \subset U, \overline{U} \cap \overline{W}=\emptyset$ and $U$ is hyperbolic.
Assume that there exists $\de \in (0, 1)$ such that for every  $f \in \Hol (\De, Y)$ with $f(0) \in V$ we have
$f(\De(0,\de)) \subset U.$ Then
$\kappa_Y (x, W)>0.$}
\vk0,2cm
\noindent
{\it Proof.} Choose a constant $c(\de)>0$ such that
$$\rho_\De (0, b) \ge c(\de) \rho_{\De (0,\de)} (0, b), \ \forall b \in \De(0, \de/2).$$
Fix an arbitrary point $y \in W$ with $\kappa_Y (x, y)<\de/2.$
Consider a holomorphic chain
$$\alpha:=\{x=x_0,x_1,\cdots, x_l=y; a_1, b_1, \cdots, a_l, b_l \in \De; f_1, \cdots, f_l \in \Hol(\De,Y)\}$$
that joins $x$ and $y$ such that
$$l(\alpha)=\rho_\De (a_1,b_1)+\cdots+\rho_\De (a_l,b_l)<\de/2.$$
It follows that $\rho_\De (a_j, b_j)<\de/2$ for every $1 \le j \le l.$
By composing with M\"oebius transformations of $\De$,
we may arrange so that
$a_1=\cdots=a_l=0,$ and hence $b_1, \cdots, b_l \in \De (0,\de/2).$ Let $k$ be the integer such that
$x_0, \cdots, x_{k-1} \in V$ but $x_k \not \in V$. By taking a refinement of $\alpha$ (see [7], p. 51), we may assume further that
$x_k \in U.$ Since $f_1 (0)=x_0, \cdots, f_k (0)=x_{k-1}$ are all in $V$, by the assumption
$f_1(\De(0, \de)), \cdots, f_k (\De(0, \de))$ are all included in $U$. Hence, the length $l(\al)$ of $\al$ may be estimated
from below as follows:
$$l(\al) \ge \sum_{i=1}^k \rho_\De (0, b_i) \ge c(\de) \sum_{i=1}^k \rho_{\De (0,\de)} (0, b_i)$$
$$\ge c(\de) \sum_{i=1}^k \kappa_U (x_{i-1}, x_i) \ge c(\de)\kappa_U (x, x_k) \ge
c(\de) \kappa_U (x, U \setminus V).$$
Here, the third inequality follows by applying the distance decreasing property to the map $f_i: \De(0, \de) \to U.$
This implies that
$$\kappa_Y (x, y) \ge c'(\de):=c(\de) \kappa_U (x, U \setminus V)>0.$$
The latter estimate follows from the fact that $\kappa_U$ defines the topology of $U$.
It follows that
$$\kappa_Y (x, W) \ge \min \{\de/2, c'(\de)\}>0.$$
Hence, we are done.
\vk.3cm
\no{\it Proof.} (of Theorem 3.1) First, let $\{x_n\}$ be a sequence in $X$ such that $x_n \to x$ in the
initial topology of $X$.
Let $U$ be a \nhd\ of $x$ which is isomorphic to some ball in a Banach space.
In particular, $U$ is hyperbolic.
Then for $n$ large enough, we have
$$\kappa_X (x_n, x) \le \kappa_U (x_n, x) \to 0, \text{as}\ n \to \infty.$$
Here the last statement follows from hyperbolicity of $U.$
Conversely, we fix $x \in X$ and a sequence $\{x_n\} \in X$ such that $\kappa_X (x_n, x) \to 0$ as $n \to \infty$.
We must show that $x_n \to x$ in the original topology of $X.$ Assume this is false, then by passing to a subsequence,
we can find an open \nhd\ $U$ of $x$ and an open \nhd\ $W$ of $\{x_n\}$ such that
$\overline{U} \cap \overline {W}=\emptyset.$ Furthermore, we can take $U$ to be hyperbolic.
Then we have $\kappa_X (x, W)=0.$ We also let $\{V_n\} \subset X$ be a sequence of open \nhd s of $x$ such that $V_n \downarrow x.$
Next, we choose a sequence $\{\de_n\}_{n \ge 0} \downarrow 0$ such that
$\de_0=1/2, \de_1=1/3$ and
$$\de_{n+1}<\min \Big \{\frac1{n}, r_n:=\de_n \prod_{j=0}^{n-1}\frac{\de_j-\de_n}{1-\de_j\de_n}\Big \}, \ \forall n \ge 1.$$
It follows that $r_{n+1}<\de_{n+1}<r_n.$ In particular, $r_n \downarrow 0.$
Using Lemma 3.2, we obtain a sequence $\{f_n\} \subset \Hol(\De, X)$ and points $a_n \in \De(0, r_n)$
such that
$$f_n (0) \in V_n, f_n(a_n) \not\in U, \ \forall n \ge 1.$$
We also set for each $n \ge 1$ the following finite Blaschke product
$$\theta_n(\lambda):=\frac{a_n}{r_n}\la \prod_{j=0}^{n-1}\dfrac{\de_j-\lambda}{1-\de_j\lambda},\  \forall \lambda\in\Delta.$$
Then $\theta_n\in \text{Hol}(\Delta, \Delta)$. Moreover, we have
$$\theta_n (0)=\theta_n (\de_j)=0, \ \forall 0 \le j \le n-1; \theta_n (\de_n)=a_n.$$
Finally, we define for each $n \ge 1$
$$g_n:=f_n\circ \theta_n\in \text{Hol}(\Delta,X).$$
Then by the above reasoning we have
$$g_n (\de_j)=g_n(0)=f_n(0), \ \forall n \ge 1, \forall 0 \le j \le n-1.$$
This implies that
$$\lim\limits_{n\to\infty}g_n(0)=\lim\limits_{n\to\infty}g_n(\de_j)=
\lim\limits_{n\to\infty}f_n(0)=x \quad\forall j\geq 0.$$
Since $X$ has the $1$-WVP and since $\de_n \downarrow 0$, we infer that $\{g_n\}_{n \ge 1}$ converges in
$\Hol(\De, X).$ In particular, there exists a small disk $\De(0,r_0)$ such that $g_n (\De(0, r_0)) \subset U$
for $n$ large enough.
This is impossible, since
$g_n (\de_n)=f_n (a_n)$ stays away from $U$ for every $n \ge 1.$
The proof is thereby completed.
\vk0.3cm
\no
{\bf Remarks.} (a) There exists a bounded Reinhardt domain $X$ in $\mathbb C^2$ that does not have 1-WVP. Indeed, let $X$ is the punctured unit polydisc in $\mathbb C^2$ i.e., $X:=\De^2 \setminus \{(0,0)\}.$
Let $\{a_k\} \in \De$ be a sequence of distinct points such that $a_k \to a \in \De \setminus \{0, 1/2, 1/3,\cdots, 1/j, \cdots\}.$
Let $f_k: \De \to X, k \ge 1$ be defined by
$$f_k (\la):=\Big (\frac{a_k-\la}{1-\la \bar a_k}, \frac{a_{k+1}-\la}{1-\la \bar a_{k+1}} \Big), \la \in \De.$$
Then $f_k \in \Hol(\De,X).$
It is also easy to check that $f_k (\la)$ converges uniformly on compact sets of $\De$ to
$$f(\la):=\Big (\frac{a-\la}{1-\la \bar a}, \frac{a-\la}{1-\la \bar a}\Big ).$$
In particular $f_k (1/j) \to f(1/j) \in X, f_k (0) \to f(0) \in X$ as $k \to \infty.$
Notice, however, that $f_k (a) \to f(a)=(0,0) \not\in X.$
Thus $X$ does not have 1-WVP.

\noindent
(b) The above domain is not pseudoconvex. In fact, we will show that every domain having WVP in $\mathbb C^n$ must be pseudoconvex. See the remark after Theorem 3.14. On the other hand,
there exists a bounded pseudoconvex domain in $\mathbb C^2$ which does not have WVP. See the remark following Proposition 4.5.
\vk0,2cm
\noindent
Our next main result provides, on the positive side, a partial converse to Theorem 3.1.
\vk0,2cm
\noindent
\textbf{Theorem 3.3.} {\it Let $X$ be a complete hyperbolic Banach analytic manifold. Then $X$ has SVP.
In particular, every compact hyperbolic manifold has SVP.}
\vk.3cm
\no
For the proof of Theorem 3.3 we need some lemmas. The first one is essentially taken from [6].
\vk.5cm
\no \textbf{Lemma 3.4.} {\it Let $X, A$ be Banach analytic manifolds and  $\{f_k\}\subset  \text{Hol}(A, X)$. Assume that $X$ is hyperbolic. Then
the sequence $\{f_k(\lambda)\}$ is a $\kappa_X$-Cauchy sequence in $X$ for every $\la \in \overline{Z_{\{f_k\}}}.$}
\vk.3cm
\no {\it Proof.}  Choose a sequence $\{\lambda_j\}\subset Z_{\{f_k\}}$ such that $\lim\limits_{j\to\infty}\lambda_j=\lambda$.
By the decreasing property of Kobayashi distance, for every $k, j, m \ge 1$ we obtain
$$\begin{aligned}\kappa_X\big(f_k(\lambda), f_m(\lambda)\big)&\leq \kappa_X\big(f_k(\lambda), f_k(\lambda_j)\big)+\kappa_X\big(f_k(\lambda_j), f_m(\lambda_j)\big)+\kappa_X\big(f_m(\lambda_j), f_m(\lambda)\big)\\
&\leq 2\kappa_{A}(\lambda_j, \lambda)+\kappa_X \big(f_k(\lambda_j), f_m(\lambda_j)\big).
\end{aligned}$$
Hence, $\{f_k(\lambda)\}$ is a $\kappa_X$-Cauchy sequence in $X$.
\vk.5cm
\no
The next lemma, a variant of Vitali's theorem for holomorphic vector-valued functions, is essentially contained in Theorem 2.1 of [1]. We include it here for the sake of completeness.

\no\textbf{Lemma 3.5.} {\it Let $E, F$ be Banach spaces and $\Omega$ be an open subset of $E.$ Let $\{f_k\}$ be a sequence in $\text{Hol}(\Om, F)$ that satisfies the following conditions:

\noindent
(i) $\{f_k\}$ is locally uniformly bounded on $\Om.$

\noindent
(ii) $Z_{\{f_k\}}$ is a set of uniqueness for $\Hol(\Om, \mathbb C)$ i.e., every holomorphic function $g: \Om \to \mathbb C$ that vanishes on $Z_{\{f_k\}}$ must be identically $0.$

Then $\{f_k\}$ converges in $\text{Hol}(\Om, F)$.}
\vk.3cm
\no{\it Proof.} Let $l^\infty (F)$ be the space of bounded sequences in $F$ equipped with the sup norm and $c(F)$ be the closed subspace of convergence sequence in $F$.
Define the map
$$f: \Om \to l^\infty (F), f(z)=(f_1 (z), \cdots, f_k (z), \cdots).$$
We split the proof into some steps.

\no
{\it Step 1.} We show that $f$ is holomorphic on $\Om.$ First, we treat the case where $n=1.$ Fix $z_0 \in \De.$ Since $\{f_k\}$ is locally uniformly bounded,
by Cauchy's inequalities we infer that $\al:= (f'_1 (z_0), \cdots, f'_k (z_0), \cdots) \in l^\infty (F).$
Now we claim that $f' (z_0)=\al.$
Using Cauchy integral formula we obtain, for $h$ small enough
$$\frac1{h} (f_k (z_0+h)-f_k (z_0))-f'_k (z_0) = \frac{h}{2\pi i} \int_{\vert z\vert=r} \frac{f_k (z)}{(z-z_0)^2(z-z_0-h)}dz,$$
where $r \in (0, 1)$ is chosen such that $\vert z_0\vert<r, \vert z_0 \vert+\vert h\vert<r.$
Since $f_k$ is locally uniformly bounded, we have
$$M:=\sup_{k \ge 1, \vert z\vert =r} \Vert f_k (z)\Vert <\infty.$$
It follows that
$$\Big \Vert \frac1{h} (f_k (z_0+h)-f_k (z_0))-f'_k (z_0) \Big \Vert \le \frac{\vert h\vert M}{(r-\vert z_0\vert)^2(r-\vert z_0\vert -\vert h\vert)}, \ \forall k \ge 1.$$
So
$$\lim_{h \to 0} \sup_{k \ge 1} \Big \Vert \frac1{h} (f_k (z_0+h)-f_k (z_0))-f'_k (z_0) \Big \Vert=0.$$
This implies $f' (z_0)=\al$ as claimed.
Thus $f$ is holomorphic on $\Om.$ For the general case,
by the above argument, $f$ is G\^ateaux holomorphic on $\Om$. Since $f$ is locally bounded, we infer that
$f$ is indeed holomorphic on $\Om,$ see Corollary II.5.5 in [2].

\no
{\it Step 2.} Let $\theta: l^\infty (F) \to l^\infty (F)/c(F)$ be the canonical projection map. We will show that
the map $f^*:= \theta \circ f$ from $\Om$ to the Banach space $\tilde F:=l^\infty (F)/c(F)$ is identically $0.$ Indeed, obviously $f^*$ is holomorphic. Moreover, by the assumption
$f^*=0$ on $Z_{\{f_k\}}$. It follows $\mu \circ f^*=0$ on $Z_{\{f_k\}}$ for every $\mu \in \tilde F',$ the dual space of $\tilde F$.
Since $\mu \circ f^*$ is a scalar holomorphic function on $\Om$ and since $Z_{\{f_k\}}$ is a set of uniqueness for $\Hol (\Om, \mathbb C)$,
we conclude that $\mu \circ f^* \equiv 0$ on $\Om.$
By Hahn-Banach's theorem, $f^*=0$ on $\Om.$

\no
{\it Step 3.} $f_k$ is convergent in $\text{Hol}(\Om, F).$ By Step 2, the sequence
$\{f_k\}$ is pointwise convergence on $\Om$.
Since $\{f_k\}$ is locally uniformly bounded on $\Om,$ we infer that the sequence
$\{f_k\}$  is equicontinuous on every compact subset of $\Om.$  Therefore, $\{f_k\}$ converges to $f$ in  $\text{Hol}(\Om, F).$
\vk0.2cm
\no
The next lemma is quite standard, it says roughly that a family of holomorphic mappings into a hyperbolic Banach analytic manifold is equicontinuous
\vk.3cm
\no\textbf{Lemma 3.6.} {\it Let $A, X$ be Banach analytic manifolds and $\{f_k\}$ be a sequence in  $\text{Hol}(A, X)$.
Assume $X$ is hyperbolic and there exists a sequence $\{\la_k\} \to \la_0 \in A$ such that $f_k (\la_k) \to x_0 \in X.$ Then the following assertions hold:

\noindent
(i) For every \nhd\ $V$ of $x_0 \in X,$
there exists an open \nhd\ $U$ of $\la_0$ in $A$ and $k_0\ge 1$ such that
$$f_k (U) \subset V,  \forall k \ge k_0.$$
\noindent
(ii) $f_k (\la_0) \to x_0$ as $k \to \infty.$}
\vk0,2cm
\no
{\it Proof.}
(i) Assume the conclusion is false. Then we may choose a sequence $\{\beta_j\} \to \la_0$ and $k_j \uparrow \infty$ such that
$$f_{k_j} (\beta_j) \not \in V\  \forall j \ge 1.$$
By the decreasing property of Kobayashi distance we get
$$\kappa_X (f_{k_j} (\beta_j), f_{k_j} (\la_{k_j})) \le \kappa_{A} (\beta_j, \la_{k_j}) \to 0 \  \text{as}\  j \to \infty.$$
It follows, using the triangle inequality, that
$$\kappa_X (f_{k_j} (\beta_j), x_0) \le \kappa_X (f_{k_j} (\beta_j), f_{k_j} (\la_{k_j})) +\kappa_X (f_{k_j} (\la_{k_j}), x_0)  \to 0 \  \text{as}\  j \to \infty.$$
This contradicts hyperbolicity of $X$. We are done.

\noindent
(ii) Applying again the triangle inequality we obtain for $k \ge 1$
the following estimates
$$\kappa_X (f_k (\la_0), x_0) \le \kappa_X (f_k (\la_k), x_0)
+\kappa_X (f_k (\la_0), f_k (\la_k))
\le \kappa_X (f_k (\la_k), x_0)+\kappa_A (\la_0, \la_k).$$
This implies that $\kappa_X (f_k (\la_0), x_0)\to 0$ as $k \to \infty$. The desired conclusion now follows from hyperbolicity of $X.$
\vk.3cm
\no
Using a standard compactness argument and Lemma 3.6 (ii) we obtain easily the following result that will be needed later on.
\vk0,2cm
\noindent
{\bf Lemma 3.7.} {\it Let $A, X$ be Banach analytic manifolds, $X$ is hyperbolic. Let $\{f_k\}$ be a sequence in
$\text{Hol}(A,X)$ which is not compactly divergent. Then there exists $\la_0 \in A$ and a subsequence $\{f_{k_j}\}$
such that $f_{k_j} (\la_0) \to x_0 \in X$ as $j \to \infty$.}
\vk0,2cm
\noindent
The following useful fact about propagation of domains on which a sequence of holomorphic maps is compactly divergent will  only be used at the end of this section.
\vk0,2cm
\noindent
{\bf Lemma 3.8.} {\it Let $A, X$ be Banach analytic manifolds and $\{\Om\}_{\al \in I}$ be a family of open subsets of $A$. Let $\{f_k\}$ be a sequence in $\text{Hol}(A,X)$ which is compactly divergent on
$\Om_\al, \forall \al \in I.$ Then $\{f_k\}$ is compactly divergent on $\Om:=\cup_{\al \in I} \Om_\al.$}
\vk0,2cm
\noindent
{\it Proof.} Suppose that $\{f_k\}$ is not compactly divergent on $\Om$. Then, there exist compact sets
$K \subset \Om, L \subset X$ and a subsequence $\{f_{k_j}\}$ such that $f_{k_j}(K) \cap L \ne\emptyset$
for every $j.$
Using compactness, we can find a sequence $\{\la_j\} \to \la_0 \in K$
such that $f_{k_j} (\la_j) \in L$ for every $j$ and $f_{k_j} (\la_j) \to x_0 \in L.$
Choose $\al_0 \in I$ such that $\la_0 \in \Om_{\al_0}.$ We may assume that the compact set
$K':=\{\la_j\} \cup \{\la_0\} \subset \Om_{\al_0}.$ Hence $f_{k_j} (K') \cap L \ne \emptyset$ for every $j$.
It follows that $\{f_k|_{\Om_{\al_0}}\}$ is {\it not} compactly divergent.
We are done.
\vk0,2cm
\noindent
Now we proceed to the proof of Theorem 3.3.
\vk.3cm
\no{\it Proof of Theorem 3.3.} Let $A$ be a connected Banach analytic manifold and
$\{f_k\} \subset \text{Hol}(A, X)$ be such that there exists some point $\la_0 \in Z^u_{f_k}$. By Lemma 3.4 and the assumption that $X$ is complete hyperbolic,  we have
$f_k (\la_0) \to x_0 \in X$ as $k \to \infty.$
Take a \nhd\ $V$ of $x_0$ in $X$ which is isomorphic to some ball in a Banach space.
Since $X$ is hyperbolic, by Lemma 3.6(i), we can find an open \nhd\ $U_0$ of $\la_0 \in A$ and $k_0 \ge 1$ such that
$$f_k (U_0) \subset V, \ \forall k\ge k_0.$$
Now we apply Lemma 3.5 to deduce that the sequence $\{f_k|_{U_0}\}$ is convergent in $\text{Hol}(U_0, X)$.
Put
$$\Omega:= \bigcup \Big \{U \subset A: \{f_k|_U\}\   \text{is convergent in}\  \Ho(U, X) \Big \}.$$
Clearly $\Omega$ is open and $U_0 \subset \Omega$ by the above proof. It suffices to show $\Om$ is closed.
Assume otherwise, then we can find $\la_1 \in \partial \Om.$
Using again Lemma 3.4 we find that $f_k (\la_1)  \to x_1 \in X$ as $k \to \infty.$ Repeating the above argument, we can find a \nhd\ $U_1$ of $\la_1$ in $A$ such that
$\{f_k|_{U_1}\}$ is convergent in $\Hol(U_1,X)$. It follows that $\lambda_1 \in U_1 \subset \Om,$ which is absurd.
Thus $\Om=A.$ The proof is complete.
\vk.3cm
\no
Our next result says that 1-WVP is in fact equivalent to WVP. We do not know if the analogous statement is true for SVP and 1-SVP.
\vk0,2cm
\noindent
\textbf{Theorem 3.9.} {\it Let $X$ be a Banach analytic manifold. If $X$ has $1$-WVP then $X$ has WVP.}
\vk.3cm
\no{\it Proof.} Let $A$ be a Banach analytic manifold and
 $\{f_k\}$ be a sequence in $\Ho(A, X)$ such that there exists
$\la_0 \in Z^u_{\{f_k\}}$ satisfying
$f_k(\la_0) \to x_0 \in X.$ Choose a \nhd\ $V$ of $x_0 \in X$ such that $V$ is isomorphic to some ball in a Banach space.
By Theorem 3.1, $X$ is hyperbolic, so we may apply Lemma 3.6(i) to find an open \nhd\ $U_0$ of $x_0$ such that $f_k (U_0) \subset V$ for $k$ large enough.
Next, by Lemma 3.5 we see that $\{f_k|_{U_0}\}$ is convergent in $\text{Hol}(U_0, X)$.

Now we set
$$\Omega:= \bigcup \Big \{U \subset A: \{f_k|_U\}\   \text{is convergent in}\  \Ho(U, X) \Big \}.$$
Clearly $\Omega$ is open and non-empty. It remains to check that $\Om$ is closed. Assume otherwise, then there exists $\la_1 \in \partial \Om$. Choose a small \nhd\  $\mathbb B$ of $\la_1$ for which we may assume to be a ball in some Banach space $E$. Pick $r>0$ and $\la_2 \in \Om$ such that
$$\la_1 \in \mathbb B(\la_2, r) \subset \Omega \cap \mathbb B \subset E.$$
Let $l$ be the  complex line  joining $\la_1$ and $\la_2.$ We may identify $\De':=l \cap \mathbb B(\la_2, r)$ with $\De.$
Consider the restriction $g_k:=f_k\big \vert_{\De'}$. By the definition of $\Om$ we see that $g_k$ is pointwise convergent on the nonempty open subset $\De' \cap  \Om$ of the disk $\De'.$
Since $X$ has $1-$WVP we infer that $g_k$ is pointwise convergent on $\De'$. In particular $f_k (\la_1)=g_k (\la_1)$ is convergent in $X$. Note that $\la_1 \in \partial \Om,$ so
$\la_1 \in Z^u_{\{f_k\}}$. Using the same reasoning as in the beginning of the proof, we see that there exists some small \nhd\  $U_1$ of $\la_1$ such that $\{f_k|_{U_1}\}$ is convergent in $\text{Hol}(U_1, X)$. Thus $\la_1 \in U_1 \subset \Om$. This contradicts the fact that $\la_1 \in \partial \Om.$
Hence $\Om=A$ and the proof is thereby completed.
\vk.3cm
\no
The next result relates weak Vitali property of a Banach analytic manifold with the usual taut property.
\vk0,2cm
\noindent
{\bf Theorem 3.10.} {\it Let $X$ be a Banach analytic manifold. Then the following statements are equivalent:

\noindent
(i) $X$ has $WVP.$

\noindent
(ii) $X$ is hyperbolic and every sequence $\{f_k\} \subset \Hol(A,X)$
where $A$ is a separable Banach analytic manifold,
contains a subsequence which is either convergent in $\Hol(A,X)$ or compactly divergent on an open dense subset of $A.$}
\vk0,2cm
\noindent
{\it Proof.} $(i) \Rightarrow (ii).$ Fix a separable Banach analytic manifold $A$ and a sequence
$\{f_k\} \subset \Hol(A,X).$ Suppose that $\{f_k\}$ contains no convergent subsequence.
Since $A$ is separable, we can choose a topological base $\{U_j\}_{j \ge 1}$ for $A.$
With no loss of generality, we can assume that each $U_j$ is a ball in some Banach space.
Fix $j \ge 1$, we claim that there exists an open non empty subset $V_j$ of $U_j$ and a subsequence $\{f_k\}_{k \in N_j}$
which is compactly divergent subsequence on $V_j.$
Indeed, assume this is false. Then we let $\{x_{j,l}\}$ be a countable dense subset of
$U_j$.
For each $l \ge 1$, we choose a small ball $\mathbb B_{j,l}:=\mathbb B(x_{j,l}, r_l) \subset U_j$ with $r_l \to 0$ as
$l \to \infty.$
Then, since $X$ is hyperbolic, we may apply by Lemma 3.7 to find,  on each ball $\mathbb B_{j,l}$ a point $y_{j,l}$ and a subsequence $\{f_k\}_{k \in N_{j, l}}$ such that $f_k (y_{j, l})$ is convergent as $k \to \infty (k \in N_{j,l})$.
Moreover, we can choose these sequences in such a way that $N_{j, l+1} \subset N_{j,l}, \forall l \ge 1$
i.e., $\{f_k\}_{k \in N_{j,l}}$ is used to construct the further subsequence $\{f_k\}_{k \in N_{j,l+1}}.$
Hence, after a diagonal process, we can build a subsequence $\{f_k\}_{k \in N_j}$ which is {\it pointwise convergence} on the countable set $A_j=\{y_{j,l}\}$ which is also {\it dense} in $U_j$.
By applying WVP of $X$ to $\{f_k|_{U_j}\}_{k \in N_j}$ we see that the sequence $\{f_k\}_{k \in N_j}$ is convergent
in $\Hol(U_j, X)$. Thus, this sequence must converge also in $\Hol(A, X)$ since $X$ has WVP.
This is absurd.
Therefore, for each $j \ge 1$, we can find an open subset $V_j$ of $U_j$ and a subsequence $\{f_k\}_{k \in N_j}$
which is compactly divergent on $V_j$. As before, we may also arrange so that $N_{j+1} \subset N_j$ for every $j \ge 1.$
Then, using one more diagonal process, we can construct a subsequence $\{f_k\}_{k \in J}$ which is compactly divergent on
{\it each} $V_j.$ Now we let
$\Om:=\cup_{j \ge 1} V_j$. Then obviously $\Om$ is open. Moreover, $\bar \Om=A$, since otherwise we would find
$j_0 \ge 1$ such that
$$V_{j_0} \subset U_{j_0} \subset A \setminus \bar \Om,$$
which is impossible. Finally, by Lemma 3.8, we conclude that $\{f_k\}_{k\in J}$ is compactly divergent on $\Om.$

\noindent
$(ii) \Rightarrow (i).$ In view of Theorem 3.9, it suffices to prove that $X$ has $1-$WVP.
For this, let $\{g_k\} \in \Hol(\De,X)$ be a sequence such that
$Z_{\{g_k\}} \cap Z^u_{\{g_k\}} \ne \emptyset.$ We have to show that $\{g_k\}$ is convergent in $\Hol(\De, X)$.
First, we claim that $\{g_k\}$ has a convergent subsequence in $\Hol (\De, X)$. Suppose otherwise, then there exist a dense open subset $\Om$ of $\Delta$ and a subsequence $\{g_{k_j}\}$ which is compactly divergent on $\Om.$
Fix $\la_0 \in Z_{\{g_k\}} \cap Z^u_{\{g_k\}}$. Then $g_{k_j} (\la_0) \to x_0 \in X.$ Since $X$ is hyperbolic, by Lemma 3.6 (i), we can choose a complete hyperbolic \nhd\ $V$ of $x_0$ in $X$ and a \nhd\ $U$ of $\la_0$ in $\Delta$ such that
$g_{k_j} (U) \subset V$ for $j$ large enough. By Theorem 3.3, $V$ has SVP. Hence $\{g_{k_j}\}$ is
convergent in $\Hol(U, V)$. This yields a contradiction to compactly divergence of
$\{g_{k_j}\}$ on the open set $U \cap \Om$ which is non-empty since $\Om$ is dense in $\Delta.$
The claim now follows.
Finally, it remains to check that two (arbitrary) accumulations points $g$ and $g'$ of $\{g_k\}$ must coincide.
For this, it suffices to note that $g=g'$ on $Z_{\{g_k\}}$. Hence, the desired conclusion now follows from the assumption that $Z_{\{g_k\}} \cap Z^u_{\{g_k\}} \ne \emptyset.$
The proof is complete.
\vk0,2cm
\noindent
{\bf Remark.} The "exceptional" set $S:= A \setminus \Om$ may depend on the sequence $\{f_k\}$ even in the case $A$ and $X$
are nice manifolds.
Indeed, let $A:=\De$ be the unit disk in $\mathbb C$ and $X$ be the unit ball of a infinite dimensional Banach space.
Let $S$ be a discrete subset of $\De$ such that $\sum_{a \in S} (1-\vert a\vert)<\infty$. We will construct a sequence
$\{f_k\} \in \Hol(\Delta, X)$ which is compactly divergent on $V:= \Delta \setminus S.$ For this, we pick sequence
$\{x_k\}$ in $X$ such that $\{x_k\}$ has no convergent subsequence.
We also let $f$ be an infinite Blaschke product associated to $S$.
Then $f \in \Hol(\De, \De)$ and $f$ vanishes exactly on $S$. Then $f_k (\la):= f(\la)x_k$ is the desired sequence.
\vk0,2cm
\noindent
The statement (ii) in the above theorem can be considerably strengthen in the special case $A=\De$ and $X$ has 1-SVP.
\vk0,2cm
\noindent
{\bf Theorem 3.11.} {\it Let $X$ be Banach analytic manifolds. Then the following assertions
are equivalent:

\noindent
(i $X$ has 1-SVP.

\noindent
(ii) Every sequence $\{f_k\} \in \Hol(\De,X)$
contains a subsequence which is either convergent in $\Hol(\De,X)$ or compactly divergent outside a discrete subset of
$\De.$}
\vk0,2cm
\noindent
{\it Proof.} $(i) \Rightarrow (ii).$ We assume that there exists no subsequence of $\{f_k\}$ which is convergent in $\Hol(\De,X).$
Let $\{r_j\}_{j \ge 1}$ be a sequence of positive numbers with $r_j\uparrow 1.$
Set $\De_j:=\De(0,r_j).$
We will prove by induction on $j$ the following statement: There exist a {\it finite} (possibly empty) set
$S_j \subset \De_j,$ an open disk $\De (0, r_j-1/j) \Subset \De'_j \Subset \De_j$
and a subsequence $\{f_k\}_{k \in N_j}$ of $\{f_k\}$
such that:

\noindent
(a) $\{f_k\}_{k \in N_j}$ is compactly divergent on $\De_j \setminus S_j;$

\noindent
(b) $N_{j+1} \subset N_j;$

\noindent
(c) $S_j \subset S_{j+1}, S_{j+1} \setminus S_j \subset \De_{j+1}\setminus \bar \De'_j.$

For $j=1,$ if the entire sequence $\{f_k\}$ is compactly divergent on $\De_1$ then we can take
$S_1:=\emptyset$ and $N_1:=\mathbb N.$ Otherwise, by Theorem 3.1 $X$ is hyperbolic, so
we use Lemma 3.7 to find
$a_{1,1} \in \De_1$ and a subsequence $\{f_k\}_{k \in N_{1,1}}$ such that $f_k (a_{1,1})$ is convergent as
$k \to \infty, k \in N_{1,1}.$ Now, if the above subsequence is compactly divergent on $\De_1 \setminus \{a_{1,1}\}$ then
we can choose $S_1:=\{a_{1,1}\}$ and $N_1:=N_{1,1}$. Otherwise, we may apply again Lemma 3.7 to get
$a_{1,2} \in \De_1\setminus \{a_{1,1}\}$ and a further subsequence $\{f_k\}_{k \in N_{1,2}}, N_{1,2}\subset N_{1,1}$
such that $f_k (a_{1,2})$ is convergent as $k \to \infty, k \in N_{1,2}.$
We claim that this process cannot be infinite. Assume on the contrary, then
we would get a sequence $\{a_{1,l}\} \subset \De_1$ of distinct points, a collection of subsequences
$\{f_k\}_{k \in N_{1, l}}, N_{1, l+1} \subset N_{1,l}$ such that
$f_k (a_{1, l})$ is convergent as $k \to \infty, k \in N_{1,l}$ for every $l \ge 1.$
Thus, using a diagonal process, we obtain a subsequence $\{f_k\}_{k \in M_1}$ such that $f_k (a_{1,l})$ is convergent for each $l \ge 1$ as $k \to \infty, k \in M_1.$
After, passing to a subsequence we may assume that $a_{1,l} \to a_1 \in \bar \De_1 \subset \De.$
Thus, using 1-SVP of $X$ we infer that the sequence $\{f_k\}_{k \in M_1}$ is convergent in $\Hol(\De,X),$ a contradiction.
Hence, the procedure described above must be finite. Thus, we can find a finite set (possibly empty)
$S_1 \subset \De_1,$ a subsequence $\{f_k\}_{k \in N_1}$ of $\{f_k\}$
such that $\{f_k\}_{k \in N_1}$ is compactly divergent on $\De_1 \setminus S_1.$

Next, suppose that there exist a finite set $S_j \subset \De_j$ and a subsequence $\{f_k\}_{k \in N_j}$
which is compactly divergent on $\De_j \setminus S_j.$ Choose a disk $\De'_j \Subset \De_j$ centered at $0$ with radius
$>r_j-1/j$ such that $S_j \subset \De'_j.$
Then by applying the preceding argument, this time, to $\{f_k\}_{k \in N_j}$ instead of
the original one $\{f_k\}$ and $\De_1$ is replaced by the annulus $\De_{j+1} \setminus \bar \De'_j,$
we obtain a subsequence $\{f_k\}_{k \in N_{j+1}}, N_{j+1} \subset N_j$ and a finite set
$S'_j \subset \De_{j+1} \setminus \bar \De'_j$ such that $\{f_k\}_{k \in N_{j+1}}$ is compactly divergent
on $\De_{j+1} \setminus (\bar \De'_j \cup S'_j).$ Then we let $S_{j+1}:=S_j \cup S'_j.$
Since
$$\De_{j+1} \setminus S_{j+1} =(\De_{j+1} \setminus (\bar \De'_j \cup S'_j)) \cup (\De_j \setminus S_j),$$
by Lemma 3.8, we infer that $\{f_k\}_{k \in N_{j+1}}$ is compactly divergent on
$\De_{j+1} \setminus S_{j+1}.$

Thus, we have proved the statement made at the beginning of the proof.
Hence, in view of (a) and (b) we may apply a diagonal process
to obtain a subsequence $\{f_k\}_{k \in I}$ which is compactly divergent on each domain $\De_j \setminus S_j, j \ge 1.$
Finally, we set $S:=\cup_{j \ge 1} S_j.$ Then, using (c) we can check that $S$ is a discrete (possibly empty) subset of $\De.$ Moreover,
since $\De \setminus S=\cup_{j \ge 1} (\De_j \setminus S_j)$, using Lemma 3.8, we deduce that $\{f_k\}_{k \in I}$ is compactly divergent on $\De \setminus S.$
The desired conclusion now follows.

\noindent
$(ii) \Rightarrow (i).$ Let $\{f_k\}$ be a sequence in $\Hol(\De,X)$ such that $Z^u_{\{f_k\}} \ne \emptyset.$
Suppose that $\{f_k\}$ contains no subsequence which is convergent in $\Hol(\De,X)$. Then, there exists a subsequence
$\{f_{k_j}\}$ which is compactly divergence outside a discrete subset $S$ of $\De.$
It follows that $Z_{\{f_k\}} \subset Z_{\{f_{k_j}\}} \subset S.$
Hence $Z^u_{\{f_k\}}=\emptyset,$ a contradiction. Thus, $\{f_k\}$ contains a convergent subsequence in $\Hol(\De,X).$
It remains to show that any two accumulations points $g$ and $g'$ of this sequence must be identical.
For this, it suffices to note that $g=g'$ on $Z_{\{f_k\}}$ and that $Z^u_{\{f_k\}} \ne \emptyset.$
\vk0,2cm
\noindent
{\bf Remarks.} (a) In [6], a Banach analytic manifold with the property described in (ii) is termed {\it weakly taut}.
Thus, Theorem 3.11 essentially generalizes (with a simpler proof) Theorem 4.1 in [6],
since the latter result is proved in the case where $X$ is a {\it finite} dimensional complex space.

\noindent
(b) It was proved in Theorem 3.4 of [6] that every complete hyperbolic Banach analytic manifold is weakly taut.
This statement also follows from our Theorem 3.3 and Theorem 3.11. Notice that our proofs does not use Zorn's lemma as
in [6].
\vk0,3cm
\noindent
Our next two results contain simple observations about inheritance of Vitali properties under inclusion.
\vk0,2cm
\noindent
\textbf{Proposition 3.12.} {\it Let $X$ be an open subset of a Banach analytic manifold $Y$. Assume that $X$ has WVP and $Y$ has SVP. Then $X$ has SVP.}
\vk.2cm
\no
{\it Proof.} Let $A$ be a Banach analytic manifold and
$\{f_k\} \subset \text{Hol}(A, X)$ be such that $Z^u_{\{f_k\}} \ne \emptyset$.
Since $Y$ has SVP we deduce that $\{f_k\}$  is convergent to $f \in \Ho(A, Y).$ Put
 $\Om:=f^{-1} (X)$. Then $\Om$ is open and non-empty since $Z_{\{f_k\}} \subset \Om.$
Notice that $\Om \supset Z_{\{f_k\}} \cap  Z^u_{\{f_k\}}$. Since $X$ has WVP,
$\{f_k\}$ is convergent in $\Ho(A,X).$ Thus, $X$ actually has SVP. We are done.
\vk0.2cm
\no
{\bf Corollary 3.13.} {\it Let $X$ be an open subset of a complete hyperbolic Banach analytic manifold $Y$.
Then $X$ has WVP if and only if $X$ has SVP. In particular, WVP and SVP are equivalent in the classes of
bounded open subsets in  Banach spaces.}
\vk0.2cm
\no
{\it Proof.} Observe that $Y$ has SVP by Theorem 3.3.
So the first assertion follows from Proposition 3.12. Finally, since every ball in a Banach space is complete hyperbolic, we get the last statement of the corollary.
\vk0.3cm
\no
This section ends up with the following result which says roughly that 1-WVP is not much weaker than SVP in the class of (finite dimensional) complex manifolds.
\vk0.3cm
\no
{\bf Theorem 3.14.} {\it Let $X$ be a complex manifold. Then the following statements are equivalent:

\no
(i)  $X$ has $1-$WVP (and hence $X$ has WVP by Theorem 3.6).

\no
(ii) $\Hol(A,X)$ is normal for every connected, locally separable Banach analytic manifold $A.$

\noindent
(iii) $X$ has SVP for source spaces $A$ having the property described in (ii).
}
\vk0.2cm
\no
{\it Proof.} $(i) \Rightarrow (ii).$ Assume that $X$ has $1-$WVP. Fix a connected, locally separable Banach analytic manifold.
Let $\{f_k\}$ be a sequence in $\Ho(A, X)$. Suppose that $\{f_k\}$ is not compactly divergent. Then, by Lemma 3.7, we can find a sequence
$\la_j \to \la_0 \in A$ and a subsequence
$f_{k_j}$ such that $f_{k_j} (\la_j) \to x_0 \in X$ as $j \to \infty.$ Let $V$ be a \nhd\ of $x_0$ which is isomorphic to some ball in an Euclidean space $\mathbb C^N.$
By Theorem 3.1, $X$ is hyperbolic, so using Lemma 3.6 (i), we can find a \nhd\ $U$ of $\la_0$ and $j_0 \ge 1$ such that
$$f_{k_j} (U) \subset V, \ \forall j \ge j_0.$$
Since $A$ is locally separable, after shrinking $U$ if necessary, we can find a countable dense subset
$Z_{\la_0}$ of $U.$
By a diagonal process, we can find a further subsequence $\{f_{k_{j_l}}\}$ which is pointwise convergence on
$Z_{\la_0}.$ It follows that $\la_0 \in Z_{\{f_{k_{j_l}}\}} \cap  Z^u_{\{f_{k_{j_l}}\}}.$
Hence, $\{f_{k_{j_l}}\}$ is convergent in $\Hol(A, X).$ We are done.

\no
$(ii) \Rightarrow (iii).$ Let $A$ be a connected, locally separable Banach analytic manifold.
Fix a sequence $\{f_k\}$ in $\Ho(A,X)$ such that
$Z^u_{\{f_k\}} \ne \emptyset.$ In particular, $\{f_k\}$ is pointwise convergence at some point of $A.$
Since $\Ho(A,X)$ is normal, we infer that $\{f_k\}$ is relatively compact in $\Ho(A,X).$
Notice that any two accumulation points of the sequence $\{f_k\}$ must be identical on $A$ in view of
the assumption that $Z^u_{\{f_k\}} \ne \emptyset.$
Therefore $\{f_k\}$ is convergent in $\Ho(A,X)$ as desired.

\no
$(iii) \Rightarrow (i)$ follows by taking $A=\De.$
\vk0.2cm
\no
{\bf Remarks.} (a) In view of the implication $(i) \Rightarrow (ii),$ we infer that every complex manifold $X$ having
WVP is necessarily taut. In particular, $X$ must be pseudoconvex at least in the case where it is a domain in
$\mathbb C^n.$ See Theorem 5.2.1 in [7].

\noindent
(b) The assumption on local separability of $A$ cannot be omitted in the implication $(i) \Rightarrow (ii).$
To see this, we consider the case where $A$ is the unit ball of $l^\infty$ and
$X=\De.$ Then, we consider the sequence of (linear) projections
$$f_k: A \to X, f_k(\la):=\la_k, \  \la=(\la_1,\cdots,\la_k,\cdots), k \ge 1.$$
Since $\{f_k\}$ contains no subsequence which is {\it pointwise} convergence on $A$
and since $\{f_k\}$ is convergent at the origin, we infer that $\Hol(A, X)$ is not normal.
\vk1cm
\no
\centerline{\bf IV.} {\bf Some classes of spaces having WVP and SVP}
\vk.3cm
\no
In this section we will investigate sufficient conditions so that a Banach analytic manifolds has Vitali properties.
For this purpose, we introduce the following terminology.

\no
{\bf Definition 4.1.} {\it An open subset $\Om$ of a Banach analytic manifold $X$ is said to have the quasi strong Vitali property (resp. quasi weak Vitali property)  if for every connected Banach analytic manifold $A$ and every sequence
$\{f_k\} \subset \Ho(A,\Omega)$ with  $Z^u_{\{f_k\}} \ne \emptyset$ (resp. with
$Z_{\{f_k\}}\cap Z^u_{\{f_k\}} \ne \emptyset$), the sequence $\{f_k\}$ is convergent in $\Hol(A,X)$.}

\noindent
These properties will be abbreviated as QSVP (resp. QWVP).
\vk0.2cm
\no
{\bf Remarks.} (a) Obviously, every open subset of a Banach analytic manifold with QSVP (resp. QWVP) also has
this property.

\noindent
(b) By Theorem 3.3, we know that every open subset of a complete hyperbolic open subset of a Banach analytic manifold has QSVP. In particular, since every ball in a Banach space is complete hyperbolic, we conclude that all open {\it bounded} subsets of a Banach space have QSVP.

\noindent
(c) Each hyperbolic relatively compact open subset $\Om$ of a complex manifold $X$
has QWVP. For this, we let $\{f_k\}$ be a sequence in $\Ho(\De,\Om)$ such that
$Z^u_{\{f_k\}} \ne\emptyset.$ Notice that, $\{f_k (z)\}$ is relatively compact in $X$ for every $z \in \De$. Furthermore, since $\Om$ is hyperbolic, the family $\{f_k\}$ is equicontinuous. By Arzela-Ascoli's theorem, $\{f_k\}$ is relatively compact in
$\Ho(\De,X).$ By the assumption that $Z^u_{\{f_k\}} \ne\emptyset,$ we deduce that two accumulation points of $\{f_k\}$
must coincide on $\De$. This implies that $\{f_k\}$ converges to some $f\in\Ho(\De, X).$
\vk0,2cm
\noindent
The first result of this section provides a class of Banach analytic manifolds having WVP.
This is a reminiscence of the well known fact that every bounded hyperconvex domain in $\mathbb C^n$ is taut
(see Corollary 5 in [9] and Proposition 5.2.2 in [7]).
\vk.2cm
\no
\textbf{Theorem 4.1}. {\it Let $X$ be a Banach analytic manifold and $\va$ be a negative \psh\ function on $X.$
Then the following assertions are equivalent:

\no
(i) $X$ is hyperbolic and for every $c<0,$ the sublevel set
$$X_c:=\{z \in X: \va(z)<c\}$$ has QWVP.

\no
(ii) $X$ has WVP.}
\vk.3cm
\no
{\it Proof.} $(ii) \Rightarrow (i)$ follows directly from Theorem 3.1.

\no
$(i) \Rightarrow (ii)$ By Theorem 3.9, it is enough to show that $X$ has $1-$WVP.
Fix a sequence $\{f_k\} \subset \Ho(\De, X)$ such that $Z_{\{f_k\}} \cap Z^u_{\{f_k\}} \ne \emptyset.$ We must show that
$\{f_k\}$ is convergent in $\Ho(\De, X).$
Choose $\la_0 \in Z_{\{f_k\}} \cap Z^u_{\{f_k\}}.$
By a reasoning as in the proof of Lemma 3.6, we can find open \nhd s $U_0 \subset \De$ of $\la_0$ and $V_0\subset X$ of $\lim_{k \to \infty} f_k (\la_0)$ and $k_0 \ge 1$ such that
$V_0$ is isomorphic to a ball in some Banach space and that
$$f_k (U_0) \subset V_0,  \forall k \ge k_0.$$
Using Lemma 3.5, we conclude that $\{f_k\}$ is convergent in $\Ho(U_0, X).$
Now we set
$$\Omega:= \bigcup \Big \{U \subset \De: \{f_k|_U\}\   \text{is convergent in}\  \Ho(U, X) \Big \}.$$
Clearly $\Omega$ is open and $U_0 \subset \Omega$ by the above proof. It suffices to show $\Om$ is closed.
Assume otherwise, then we can find $\la_1 \in \partial \Om.$ Notice that $\{f_k\}$ converges to $f \in \Ho(\Om, X).$
Now we set
$$\psi (z):= \sup_{k \ge k_0} (\va \circ f_k) (z), \ \forall z \in \De.$$
By the assumption on $\va$ we infer that the upper regularization $\psi^* \le 0$ and is subharmonic on $\De$.
Furthermore, by the choice of $U_0$ we have
$$\sup_{U_0} \psi^* \le \sup_{V_0} \va <0.$$
Thus the maximum principle yields
$\psi^*<0$ entirely on $\De.$ In particular $\psi^* (\la_1)<0$. Fix $c \in (\psi^* (\la_1), 0).$ Choose an open disk $U_1$ around $\la_1$ such that
$\sup_{U_1} \psi^* <c.$
It follows that
$$\sup_{U_1} \va \circ f_k <c, \  \forall k \ge k_0.$$
Therefore $f_k$ map $U_1$ into $X_c$ for every $k \ge k_0.$ Notice that
$$\emptyset \ne U_1 \cap \Om \subset U_1 \cap Z_{\{f_k\}} \cap Z^u_{\{f_k\}}.$$
Since $X_c$ has QWVP, we deduce that the sequence $\{f_k\}$ is convergent in $\Ho(U_1, X).$
Thus $\la_1 \in \Om$. This contradicts our choice that $\la_1 \in \partial \Om.$ The proof is therefore completed.
\vk.3cm
\noindent
{\bf Corollary 4.2.} {\it Let $X$ be a hyperbolic complex manifold. Assume that there exists a negative exhaustion function
$\va$ for $X$ i.e, $X_c:= \{z \in X: \va(z)<c\}$ is relatively compact in $X$ for every $c<0.$
Then $X$ has WVP.}
\vk0,2cm
\noindent
{\it Proof.} By the remark (c) following Definition 4.1 and the assumption on hyperbolicity of $X$, we see that all sublevel sets $X_c$ has QWVP. The desired conclusion now follows immediately from Theorem 4.1.
\vskip0,3cm
\no
The theorem below gives a sufficient condition for Vitali properties for open subsets of Banach analytic manifolds.
\vk.3cm
\no{\bf Theorem 4.3}. {\it Let $Y$ be an open subset of a Banach analytic manifold $X$. Assume that $X$ has WVP and that there exists a negative \psh\ function $\va$ on $Y$ such that
$\lim_{z \to \xi} \va(z)=0$ for every $\xi \in \partial Y.$
Then $Y$ has WVP.}
\vk.3cm
\no{\it Proof.}
According to Theorem 3.9,  it suffices to show that $Y$ has $1-$WVP. Fix a sequence
$\{f_k\}_{k\geq 1}\subset \text{Hol}(\Delta, Y)$ such that $Z_{\{f_k\}} \cap Z^{u}_{\{f_k\}} \ne \emptyset.$
Then $\{f_k\}_{k\geq 1}\subset \text{Hol}(\Delta, X)$. Thus $\{f_k\}$ is convergent to $f \in  \text{Hol}(\Delta, X)$, since $X$ has WVP.
We claim that $f(\De)\subset Y$. To this end, we set
$$u(z):=\limsup_{k \to \infty} (\varphi\circ f_k)(z),\;\; z\in \Delta.$$
Fix $x_0 \in \De \cap  Z_{\{f_k\}}.$
Then the following statements are true:

\no
(a) $u(z) \le (\varphi\circ f)(z) \  \forall z\in\Delta':=f^{-1}(Y)\cap \Delta$. This follows from upper semicontinuity of $\va$.
Notice $\Delta'$ also that is open and $x_0\in \Delta'$.

\no
(b) $u$ satisfies the sub-mean value inequality, i.e.,
$$u(z_0)\leq \dfrac1{2\pi}\int_0^{2\pi}u(z_0+re^{i\theta})d\theta,\;\; \forall z_0\in \Delta,\; \forall r>0\; \text{small enough}.$$
This is an easy consequence of Fatou's lemma.

\no The problem is to show $\De' =\De$. Assume otherwise, then $F:=\Delta\setminus \Delta'\ne\emptyset$.
Since $\varphi\circ f$ is subharmonic on $\Delta'$ and since $(\va \circ f) (x_0)<0$, by the maximum principle  we obtain
$$\varphi\circ f<0 \ \text{on}\  \Delta'.$$
Combining with (a) we get
$$u(z)<0\quad \forall z\in\Delta'.$$
Next, we pick $x_1\in \De \cap \partial F$. Then $x_1\in F$, since $F$ is closed. Hence $f(x_1) \in \partial Y.$ This implies that
$$u(x_1)=\lim_{\xi \to f(x_1)} \va(\xi)= 0.$$
\no
Choose $r>0$ so small such that the closed disk $\bar \De (x_1, r)$ is included in $\De.$ Thus
$$\partial \De(x_1, r) \cap \De' \ne \emptyset.$$
It follows that there exist $\theta_0 \in (0, 2\pi)$ and $\de>0$ such that
$$x_1+re^{i\theta} \in \De'  \  \forall \theta \in (\theta_0-\de, \theta_0+\de).$$
By (b) we obtain
$$0=u(x_1)\leq  \dfrac1{2\pi}\int_0^{2\pi}u(x_1+re^{i\theta})d\theta \le  \dfrac1{2\pi}\int_{\theta_0-\de}^{\theta_0+\de}u(x_1+re^{i\theta})d\theta<0.$$
The last inequality follows from $u<0$ on $\Delta'.$
This is absurd. Thus $\Delta'=\Delta$. Hence $f(\Delta)\subset Y$. The proof is complete.
\vk.3cm
\no
Now we consider conditions that imply SVP of open subsets of a Banach space. This result will be used to characterize
SVP of balanced domain in Banach spaces.
\vk0,2cm
\noindent
{\bf Theorem 4.4.} {\it Let $\Om$ be an open subset of a Banach  space $X$. Assume that there exists a negative \psh\ function $\va$ on $\Om$ such that
$$\lim_{z \to \xi} \va (z)=0, \  \forall \xi \in \partial \Om.$$
\no
Then the following assertions are equivalent:

\no
(i) $\Om$ is hyperbolic and for every $c<0$ the open set $\Om_c:=\{z \in \Om: \va(z)<c\}$ has QSVP.

\no
(ii) $\Om$ has SVP.}
\vk.3cm
\no
{\bf Proof.} $(ii) \Rightarrow (i)$ follows from Theorem 3.1.

\no
$(i) \Rightarrow (ii).$ Fix a connected Banach analytic manifold $A.$
Let $\{f_k\}$ be a sequence in $\Ho(A,X)$ such that $Z^u_{\{f_k\}} \ne \emptyset.$
Fix $\la_0 \in Z_{\{f_k\}}$. Then we have
$f_k (\la_0) \to x_0 \in X.$ By upper-semicontinuity of $\va$, we
can choose a \nhd\ $V$ of $x_0 \in X$ such that $\sup_V \va<0.$
Using hyperbolicity of $X$, by Lemma 3.6(i), we can find a \nhd\ $U$ of $\la_0$ in $A$ and $k_0 \ge 1$ such that
$$f_k (U) \subset V \  \forall k \ge k_0.$$
Set
$$\psi (z):= \sup_{k \ge k_0} (\va \circ f_k) (z), \ \forall z \in A.$$
Since $\va$ is negative and \psh\ on $A$, we infer that the upper regularization $\psi^*$ is also plurisubharmonic on $A$ and $\le 0$ there.
Moreover, by the choice of $U$ and $V$ we also have
$$\psi^* (\la_0) \le \sup_V \va<0.$$
Hence  the maximum principle yields
$\psi^*<0$ entirely on $A.$ Next, we choose an arbitrary point $\la_1 \in Z^u_{\{f_k\}}.$
We claim that $f_k$ is uniformly bounded on a small open \nhd\ of $\la_1.$
To see this, choose $c\in \mathbb R$ such that  $\psi^* (\la_1)<c<0.$
Then, there exists a open \nhd\ $W$ of $\la_1$ such that $\sup_{W} \psi^*<c.$
It follows that $f_k (W) \subset \Om_c$ for every $k \ge k_0.$
Since $\Om_c$ has QSVP, we can find an
open subset $Y$ of $X$ that contains $\Om_c$ such that $Y$ has SVP. Hence, the sequence $\{f_k\big\vert_{W}\}$ is convergent in $\Ho(W, Y)$.
Thus $f_k$ must be locally uniformly bounded near $\la_1$ for $k \ge k_0$.
The claim is proved. It means that we can find an open \nhd\ $\mathbb B$ of $\la_1$ on which $f_k$ is uniformly bounded for $k \ge k_0.$
Now, we use  Lemma 3.5 to conclude that $\{f_k\}$ is convergent to $f$ in $\text{Hol}(\mathbb B, X)$.
We claim that $f(\mathbb B) \subset \Om.$ Assume otherwise, then there exists $x_0 \in \mathbb B$ such that $f(x_0) \in \partial \Om.$
Since $\psi$ tends to $0$ at $f(x_0)$ we must have $\psi (x_0)=0$. This is absurd.
Thus $f(\mathbb B) \subset \Om$ as desired. Now we let
$$A':= \bigcup \Big \{U \subset A: \{f_k|_U\}\   \text{is convergent in}\  \Ho(U, \Om) \Big \}.$$
Clearly $A'$ is open, by the above reasoning $A'$ is also non-empty. It remains to check $A'$ is closed. Assume otherwise, then there exists $\la_2 \in \partial A'$.
Then $\la_2 \in  Z^u_{\{f_k\}}.$ Repeating the preceding argument, we see that there exists a small ball  $\mathbb B'$ around $\la_2$ such that $f_k$ is uniformly convergent in
$\text{Hol}(\mathbb B', \Om)$. Hence $\la_2 \in A'.$ This is impossible. The proof is therefore complete.
\vk0.3cm
\noindent
We now discuss Vitali properties of special classes of Banach analytic manifolds in the rest of this section.
The first objects to consider are Hartogs domains over Banach analytic manifolds.

Recall that, given a Banach analytic manifold $X$ and an upper semicontinuous function $\va: X \to [-\infty, \infty),$
the Hartogs domain $\Om_\va (X)$ is defined as
$$\Om_\va (X):=\{(z, w) \in X \times \mathbb C: \vert w\vert<e^{-\va(z)}\}.$$
\noindent
The next result relates Vitali properties of a Hartogs domain and those of its base and radii of fibers.

\noindent
{\bf Proposition 4.5.} {\it The Hartogs domain $\Omega_\va (X)$ has SVP (resp. WVP) if and only if $\va$ is continuous \psh\ on $X$ and $X$ has SVP (resp. WVP).}
\vk0.2cm
\noindent
{\bf Remark.} Thus, if $X=\De$ and $\va$ is bounded subharmonic but not continuous on $\Delta$ then $\Om_\va (\De)$
is a bounded pseudoconvex domain in $\mathbb C^2$ without having WVP.
\vk0,2cm
\noindent
{\bf Proof.} We only give the proof for the SVP case, the other case is similar and somewhat simpler.

\noindent
$(\Rightarrow).$ Suppose that $\Omega_\va (X)$ has SVP. First we check that $X$ has $SVP.$
For this, let $A$ be a connected Banach analytic manifold and
$\{f_k\} \in \Ho(A, X)$ be such that $Z^u_{\{f_k\}} \ne \emptyset.$
Set $f'_k:=(f_k, 0).$ It is then clear that $\{f'_k\} \in \Ho(A,\Om_\va (X))$. Moreover,
$Z^u_{\{f'_k\}} \ne \emptyset.$ It follows, using SVP of $\Om_\va (X)$ that $\{f'_k\}$ converges in $\Hol(A,\Om_\va (X))$. Thus, so does $\{f_k\}$. Hence $X$ has SVP.

Now we prove continuity of $\va$. Assume that $\va$ is discontinuous at $x^* \in X.$
Then, since $\va$ is upper semicontinuous, we can find a sequence $x_k \to x^*$ and $s \in \mathbb R$ such that
$$\va (x_k) \le s<\va (x^*), \ \forall k.$$
Next, we set $r:=e^{-s}$ and define a sequence $\{g_k\} \subset \Ho(\De,\Om_\va (X))$ by
$$g_k(\la):=(x_k, r\la), \ \forall \la \in \De.$$
If $\vert \la\vert<\de:=r^{-1}e^{-\va(x^*)}$ then $\vert \la\vert r<e^{-\va(x^*)}.$
Hence,
$$g_k (\la) \to g(\la):= (x^*, r\la) \in \Om_\va (X),$$
for $\la \in \de\De.$ In particular, $0 \in Z^u_{\{f_k\}}.$
Since $\Om_\va (X)$ has SVP, the sequence $\{g_k\}$ must converge to
$\tilde g \in\Ho(\De, \Om_\va(X)).$ By the above reasoning, $f$ agrees with
$\tilde f$ on $\de\De.$
By uniqueness property of holomorphic maps from $\De$ to $X \times \mathbb C$, we infer that $\tilde g=(x_0, r\la)$ for all $\la \in \De.$
Hence $(x^*, r\la) \in \Om_\va (X)$ for all $\la\in \De.$ This yields a contradiction to the choice of $r.$
Thus $\va$ is continuous on $X.$
It remains to prove that $\va$ is plurisubharmonic on $X$. To see this, it is enough
to show $\va$ is plurisubharmonic on every open set $U$
which is isomorphic to an open subset of a Banach space. Fix such an open set $U$ and
let $\theta: \De \to U$ be an arbitrary holomorphic map. It suffices to show that the continuous function
$u:=\va \circ \theta:\De \to \mathbb R$ is subharmonic. Assume otherwise, then we can find a closed disk $\De' \subset \De$, a holomorphic polynomial $p$ in $\mathbb C$ such that
$$u \le \Re p \  \text{on}\ \partial \De' \ \text{whereas}\ \varepsilon:=\sup_{x \in \De'} (u(x)-\Re p(x))>0.$$
For $k \ge 1,$ we define
$$h_k(\la):= (\theta(\la), e^{-p(\la)-\varepsilon-\frac1{j}}) \ \forall \la \in \De'.$$
By the choice of $\varepsilon$ we can check that
$h_k \in \Ho(\De',\Om_\va (X)).$ Furthermore, we also note
$$h_k(\la) \to h(\la):= (\theta(\la), e^{-p(\la)-\varepsilon}), \ \forall \la \in \De'.$$
Note that $h \in \Ho(\De', X \times \mathbb{C}).$
Now we choose $\al \in (0,1)$ such that
$\Re p(\la)+\varepsilon>u(\la)$ if $\al \in V_\al:=\De \setminus \De (0,\al).$
It follows that $h(V_\al) \subset \Om_\va (X)$.
Since $\Om_\va (X)$ has SVP we deduce that $h_k$ converges to
$\tilde h \in \Ho(\De', \Om_\va (X)).$ Since $h=\tilde h$ on $V_\al$,
using again uniqueness property of holomorphic maps from $\De'$ to $X \times \mathbb C$ we obtain $h=\tilde h$
on $\De'.$ This implies that
$\Re p(\la)+\varepsilon>u(\la)$ for {\it all} $\la \in \De'$. This contradiction to the choice of $\varepsilon$ proves plurisubharmonicity of $\va$ on $X.$

\noindent
$(\Leftarrow).$ Assume that $\va$ is continuous \psh\ on $X$ and $X$ has SVP. Fix a connected Banach analytic manifold $A$ and a sequence $\{f_k\} \in \Ho(A, \Om_\va (X))$ satisfying $Z^u_{\{f_k\}} \ne \emptyset.$
We write $f_k=(g_k,h_k),$ where $g_k \in \Hol(A, X)$ and $h_k \in \Hol (A,\mathbb C).$
Then we have
$$Z^u_{\{g_k\}} \ne \emptyset,Z^u_{\{h_k\}} \ne \emptyset.$$
Since $X$ has SVP, we deduce that $\{g_k\}$ converges
to $g \in \Hol(A,X).$ Notice also that
$$\vert h_k(\la)\vert<e^{-\va(g_k(\la))}, \ \forall \la \in A, \forall k \ge 1.$$
It follows that the sequence $\{h_k\}$ is uniformly bounded on compact sets of $A.$
By Lemma 3.5, we infer that $\{h_k\}$ is convergent to $h \in \Hol(A, \mathbb C).$
This implies that
$$\vert h(\la)\vert\le e^{-\va (g(\la)}, \forall \la \in A.$$
Rewriting the above inequality as
$$f(\la):=\log \vert h(\la)\vert+\va (g(\la)) \le 0, \ \forall \la \in A.$$
Since $g, h$ are holomorphic functions on $A$ and since $\va$ is \psh\ on $X$ we infer that $f$ is \psh\ on $A$.
Moreover, $f<0$ on the non-empty set $Z_{\{f_k\}}.$ It follows, using the maximum principle, that
$f(\la)<0$ for every $\la \in A.$ Therefore
$$\vert h(\la)\vert < e^{-\va (g(\la)},\ \forall \la \in A.$$
Thus $\{f_k\}$ converges to $(g, h)\in \Hol(A,\Om_\va (X)).$ Hence $\Om_\va (X)$ has
SVP as desired.
\vk0.3cm
\noindent
The next result deals with Vitali properties of balanced domains in Banach space.
Recall that a domain $\Om$ in a Banach space $E$ is said to be balanced if $x \in \Om$ then $\la x$ in $\Om$ for every
$\la \in \De.$ In particular $0 \in \Om.$ For a balanced domain $\Om$, the gauge (or Minkowski) functional of $\Om$ is defined as
$$h_\Om (x):= \inf\{\la>0: x \in \la \Omega\}, \ x \in E.$$
It is clear that $h_\Om$ is {\it homogeneous} i.e., $h(\la x)=\vert \la\vert h(x)$ and, since $\Om$ is a domain,
$h_\Om$ is upper semicontinuous and
$$\Om=\{x \in E: h_\Om (x)<1\}.$$
\noindent
We are now able to formulate the final result of this section.
\vk0,2cm
\noindent
{\bf Proposition 4.6.} {\it Let $\Om$ be a balanced domain in a Banach space $E$. Then the following statements are equivalent:

\noindent
(i) $\Om$ has WVP.

\noindent
(ii) $\Om$ is bounded and $h_\Om$ is continuous on $E$ and plurisubharmonic on $\Om.$

\n
(iii) $\Om$ has SVP.}
\vk0.2cm
\noindent
{\bf Proof.} $(i) \Rightarrow (ii).$ If $\Om$ has WVP then by Theorem 3.1, $\Om$ is hyperbolic.
Thus using the same argument as in the beginning of the proof of Theorem 6.1 in [6] we conclude that $\Om$ is bounded.
Next, we show that $\log h_\Om$ is \psh\ on $\Om.$
For this, since $h_\Om$ is upper semicontinuous on $\Om,$
it suffices to show that for every  choice $a, b \in E$ the function
$\la \mapsto u(\la):=\log h_\Om (a+\la b)$ is subharmonic on the open set $\Om_{a, b} \subset \mathbb C$ where it defines.
If this is not the case then we can find a closed disk $\De' \subset \Om_{a, b}$ and a
holomorphic polynomial $p$ in $\mathbb C$ such that
$$u \le \Re p \  \text{on}\ \partial \De' \ \text{whereas}\ \varepsilon:=\sup_{x \in \De'} (u(x)-\Re p(x))>0.$$
For $k \ge 1,$ we define
$$f_k(\la):= \frac{a+\la b}{e^{p(\la)+\ve+1/k}} \ \forall \la \in \De'.$$
By the choice of $\varepsilon$ we can check that
$f_k \in \Ho(\De',\Om).$ Furthermore, we also note
$$f_k (\la) \to f(\la):=\frac{a+\la b}{e^{p(\la)+\ve+1/k}} \ \forall \la \in \De'.$$
By the choice of $p$, we see that there exists an open \nhd\ $V$ of $\partial \De'$ such that
$f(V) \subset \Om.$ Since $\Om$ has WVP we deduce that $f_k$ converges to
$\tilde f \in \Ho(\De', \Om).$ Since $f=\tilde f$ on $V$,
using again uniqueness property of holomorphic maps from $\De'$ to $E$ we obtain $f=\tilde f$
on $\De'.$ This implies that
$\Re p(\la)+\varepsilon>u(\la)$ for {\it all} $\la \in \De'$. This contradiction to the choice of $\varepsilon$ proves plurisubharmonicity of $\log h_\Om$ on $\Om_{a, b}.$ Thus $\log h_\Om$ and hence $h_\Om$ is \psh\ on $\Om.$
It remains to check continuity of $h_\Om$ on $E$. Suppose $h_\Om$ is discontinuous at $x^* \in E$. Then, since $h_\Om$ is upper semicontinuous at $x^*$, we can find a sequence $\{x_k\} \to x^*$ and $s>0$ such that
$$h_\Om (x_k)<s<h_\Om (x^*), \ \forall k \ge 1.$$
For $k \ge 1$, we define  $$g_k (\la):=\frac{\la}{s}x_k, \ \forall \la \in \De.$$
By the choice of $s$, we have $g_k \in \Ho(\De, \Omega)$. Moreover,
$$g_k (\la)\to g(\la):= \frac{\la}{s}x^*,  \text{as}\ k \to \infty, \forall \la \in \De.$$
Since $\Omega$ contains a \nhd\ of $0,$ there exists $\de>0$ such that $g(\De(0,\de)) \subset X.$
It follows that $\De (0,\de)\subset Z_{\{g_k\}}$. Since $\Omega$ has WVP, we infer that $g(\De)$ must be included in $\Omega$.
Hence, $h_\Om (x^*) \le s$. This contradiction proves continuity (on $E$) of $h_\Om$.

\noindent
$(ii) \Rightarrow (i).$ Suppose that $\Om$ is bounded and $h_\Om$ is continuous plurisubharmonic on $\Om.$
Let $\va:= h_\Om-1$. Then $\va$ is negative \psh\ on $\Om.$ Moreover, fix $\xi \in \partial \Om$, since $h_\Om$ is continuous at $\xi$ we infer that
$\lim_{z \to \xi} \va(z)=0.$
Notice also that, being a bounded domain in a Banach space, $\Om$ has QSVP. Therefore, we may apply Theorem 4.4 to reach that $\Om$ has WVP.

\noindent
$(ii) \Rightarrow (iii).$ By the above implication $\Om$ has WVP. In view of Corollary 3.13, $\Om$ has SVP.

\noindent
$(iii) \Rightarrow (i)$ is trivial.

The proof is thereby completed.
\vk0,5cm
\noindent
{\bf V. Open questions}
\vk0,2cm
\noindent
Before leaving this paper, we wish to point out a few questions that are left open by our methods.

\noindent
1. Is there a Banach analytic manifold with WVP but without SVP? We conjecture that there exists such a Banach analytic manifold.

\noindent
2. Is there any analogue of Theorem 3.11 in the case where $X$ has SVP i.e., the sequence $\{f_k\}$ is completely divergent
outside a set which is locally contained in {\it an analytic hypersurface}?

\noindent
3. Using Proposition 4.5, we see that the Hartogs domain $\Om_\va (\De)$ is {\it unbounded} and has SVP if $\va$ is continuous, subharmonic and satisfies $\inf_\De \va=-\infty.$
Is there any substantial class of unbounded domains (in Banach spaces) having WVP and SVP?
More precisely, can we describe WVP and WVP of an unbounded domain in terms of the existence of peak plurisubharmonic functions at finite and {\it infinite} boundary points?
\vk.5cm

\vk.3cm

\no Department of Mathematics \& Informatics, Hanoi National University of Education, 136 Xuan Thuy Street, Hanoi, Vietnam.

\no Email-address: dieu\_vn@yahoo.com (Nguyen Quang Dieu).

\no Email-address: nvkhiemdhsp@gmail.com (Nguyen Van Khiem).


\begin{thebibliography}{1}

\bibitem{ArenNi} W. Arendt  and N. Nikolski, \emph {Vector-valued holomorphic functions revisited}, Math. Zeit. {\bf 234}
(2000), 777-805.

\bibitem{Din} S. Dineen, \emph {The Schwarz lemma}, Clarendon Press, 1989.

\bibitem{FraVesen} T. Franzoni and E. Vesentini, \emph {Holomorphic Mappings and Invariant Distances,} North-Holland
Mathematical Studies vol. 40, North-Holland, Amsterdam-New York-Oxford, 1980.

\bibitem{HaiBan}  L. M. Hai and P. K. Ban, \emph {On the tautness and locally weak tautness of domain in a Banach space,} Acta Math. Vietnamica,  \textbf{28} (2003), No1, 39-50.

\bibitem{HaiKhueTrang}  L. M. Hai, N. V. Khue and P. N. T. Trang, \emph {Normality of a family of Banach-valued holomorphic maps,} Acta Math. Vietnamica,  \textbf{29} (2004), No3, 251-257.


\bibitem{HaiQuangVyHung} L. M. Hai, T. T. Quang, D. T. Vy and L. T. Hung,  \emph {Some classes of Banach analytic spaces},  Mathematical Proceedings of the Royal Irish Academy, \textbf{116}A (2016), 1-17.

\bibitem{Kobayshi} S. Kobayashi, \emph {Hyperbolic Complex Spaces}, Springer 1998.

\bibitem{Pa} R. Palais, \emph {Homotopy theory of infinite dimensional manifolds,} Topology, \textbf{5} (1966), 1-16.

\bibitem{Sib} N. Sibony, \emph{A class of hyperbolic manifolds}, In: "Recent Developments in Several Complex
Variables", J. E. Fornaess (ed.), Ann. Math. Studies 100 (1981), 347-372.
\end{thebibliography}
\end{document}